\documentclass[12pt]{amsart}
\usepackage{amscd}
\usepackage{verbatim}
\usepackage{amssymb, amsmath, amsthm, amscd,ifthen}
\usepackage[dvips]{graphics}
\usepackage[cp866]{inputenc}
\usepackage{graphicx}
\usepackage{epsfig}

\usepackage{amsmath,amssymb,amscd,}
\usepackage[all]{xy}

\usepackage[dvips]{graphics}
\usepackage[cp866]{inputenc}
\usepackage{graphicx}
\usepackage{epsfig}
\usepackage[hang]{subfigure}

\theoremstyle{plain}
\newtheorem{theorem}{Theorem}
\newtheorem{lemma}{Lemma}
\newtheorem{corollary}{Corollary}
\newtheorem{proposition}{Proposition}

\theoremstyle{definition}
\newtheorem{definition}{Definition}
\newtheorem{example}{Example}

\theoremstyle{plain}
\newtoks\thehProclaim
\newtheorem*{Proclaim}{\the\thehProclaim}

\theoremstyle{definition}
\newtoks{\thehRemark}
\newtheorem*{Remark}{\the\thehRemark}

\renewcommand{\leq}{\leqslant}

\begin{document}

\title[Oriented area is a perfect Morse function]{Oriented area is a perfect Morse function}

\author{Gaiane Panina}

\address{Saint Petersburg State University,
\newline 7/9 Universitetskaya nab, St. Petersburg, 199034 Russia}

\email{gaiane-panina@rambler.ru}



\begin{abstract}

 We show that an appropriate generalization of the \textit{oriented area function}
  is a perfect Morse function on  the space of  three-dimensional configurations of
an equilateral polygonal linkage with odd number
of edges.
Therefore  cyclic equilateral polygons (which appear as  Morse points) are interpreted as
independent generators of the homology groups of the (decorated)
configuration space.
\end{abstract}

\date{11/July/2016}

\maketitle

\section{Introduction}

A Morse function on a smooth manifold is called perfect if the number of critical points equals the sum of Betti numbers.
Not every manifold possesses a perfect Morse function. Homological spheres (that are not spheres) do not possess it;
manifolds with torsions in homologies do not possess it, etc. On the other hand, the celebrated Millnor-Smale cancelation technique for critical points with neighbor indices (or, equivalently, cancelation of handles)
provides a series of existence-type theorems \cite{Fomenko,Milnor} which are the key tool in Smale's proof of generalized Poincare conjecture \cite{Smale}.

In the paper we focus on one particular example of a perfect Morse function and discuss some related problems.
Namely, we bound ourselves by the space  of configuration of the equilateral polygonal linkage with odd number $n=2k+1$ of edges. As the ambient space it makes sense to take either $\mathbb{R}^2$ or $\mathbb{R}^3$,
which gives us the spaces $M_2(n)$ and $M_3(n)$.
In bigger  dimension the configuration space is not a manifold. The number $n$ is chosen to be odd by the same reason: for even $n$, the configuration space of the equilateral polygonal linkage has singular points.

We are interested in finding a "natural" perfect Morse function, that is, a function that has a transparent physical or geometrical meaning.
The first candidate of a "natural" Morse function on  $M_2(n)$ was the oriented area function $A$. Indeed,  it a Morse function with
 an easy description of its critical points  as cyclic polygons (that is, polygons with a superscribed circle), and a simple formula for
the Morse index of a critical point \cite{khipan}. However, for $M_2(n)$
 $A$ is not  perfect. In particular, for the
equilateral pentagonal linkage it has one extra local maximum
(except for the global maximum) and one extra local minimum, see
Example \ref{Expentagon}.
  For the equilateral heptagonal linkage the number
of Morse points greatly exceeds the sum of Betti numbers of the
configuration space.

To build up a perfect Morse function,  we take the  space $M_3(n)$ and \textit{decorate} it.
The decorated space $\widetilde{M}_3(n)$ is well adjusted for
 an appropriate
generalization $S$ of the area function $A$.
Its critical points (loosely speaking) are again cyclic polygons.
Surprisingly, the function $S$ is a perfect Morse function.
 As a
direct corollary we interpret cyclic equilateral polygons as
independent generators of  the homology groups of the configuration
space  $\widetilde{M}_3(n)$.

\section{Preliminaries and notation}\label{section_preliminaries}
For an odd $n=2k+1$
an equilateral \textit{polygonal  $n$-linkage} should be interpreted as a collection of rigid
bars of lengths $1$ joined consecutively by revolving joints in a
chain.

 \textit{A configuration} of the polygonal  $n$-linkage in the
Euclidean space $ \mathbb{R}^d$, $d=2, \ 3$ is a sequence of points
$R=(p_1,\dots,p_{n+1}), \ p_i \in \mathbb{R}^d$ with
$1=|p_i,p_{i+1}|$ and $1=|p_n,p_{1}|$ modulo the action of
orientation preserving isometries of the space  $\mathbb{R}^d$.
We also call $P$ \ a \textit{polygon}. A configuration carries
a natural orientation which we indicate in figures by an arrow.

The set $M_d(n)$ of all  configurations up to an orientation-preserving isometry of the ambient space is \textit{the moduli
space, or the configuration space of the polygonal linkage }$L$.

For $d=2,3$ the space $M_d(n)$ is a smooth manifold.

 We explain below in this paragraph what is known about planar
configurations and the signed area function as the Morse function on
the configuration space.

\begin{definition} \label{Dfn_area} The \textit{signed area} of a polygon $P\in M_2(n)$ with the vertices \newline $p_i = (x_i,
y_i)$  is defined by
$$2A(P) = (x_1y_2 - x_2y_1) + \ldots + (x_ny_1 - x_1y_n).$$
\end{definition}

\begin{definition}
    A
polygon  $P$  is  \textit{cyclic} if all its vertices $p_i$
lie on a circle.

\end{definition}

 Cyclic polygons  arise here  as critical points of the signed area:
 a polygon $P$ is a critical point of the
signed area function $A$  iff $P$ is a cyclic configuration.
Their Morse indices are known \cite{khipan,zhu}.

\medskip

\begin{example} \label{Expentagon} \cite{panzh} The equilateral pentagonal linkage  has 14 cyclic
configurations  indicated in  Fig.  4.

(1). The convex regular pentagon and its mirror image  are the global maximum and
minimum of the signed area $A$. Their Morse indices are $2$ and $0$
respectively.

(2). The starlike configurations are a local maximum and a local
minimum of $A$.

(3).  There are ten  configurations that have three consecutive
edges aligned. Their  Morse indices equal $1$.

\end{example}

\begin{figure}[h]\label{pentagon}
\centering
\includegraphics[width=5 cm]{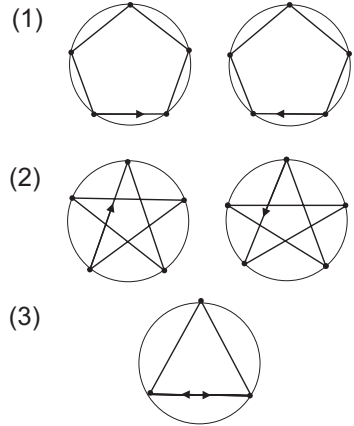}
\caption{Cyclic configurations of the equilateral pentagonal linkage
}
\end{figure}



\section{The decorated configuration space $\widetilde{\mathcal{M}}_3(n)$ and the area function $S$}\label{Section_decorated}

\begin{definition}
The \textit{decorated configuration space} is defined as the set of pairs
$$\widetilde{M}_3(n)=\{(P,\xi)| P \ \hbox {is a polygon in $\mathbb{R}^3$
 with the sidelengths}\  1; \ \xi \in S^2,
\}$$  factorized by the diagonal action of the orientation
preserving isometries of $\mathbb{R}^3$.

Here $S^2 \in \mathbb{R}^3$ is the unit sphere centered at the
origin $O$.
\end{definition}

\begin{lemma}\begin{enumerate}
               \item The space $\widetilde{M}_3(n)$ is an orientable fibration over $M_3(n)$ whose
fiber is $S^2$.

               \item The Euler class of the fibration equals zero.
             \end{enumerate}

\end{lemma}
Proof. (1) The set of all polygons with fixed sidelengths (before
factorization by isometries) is known to be orientable. Therefore the set of the
pairs (a polygon, a vector) is also orientable as a trivial
fibration. Since we take a factor by the action
of orientation preserving isometries, the result is also orientable.

(2)  $s(P):=\frac{\overrightarrow{p_1p_2}}{|p_1p_2|}$ defines  a non-zero section. \qed

The Gysin sequence implies:
\begin{corollary}\label{Cor_Gisin}
\begin{enumerate}
  \item We have the following short exact sequence:
$$0\rightarrow H^m(M(n)) \rightarrow H^m(\widetilde{M}(n)) \rightarrow H^{m-2}(M(n)) \rightarrow 0.$$
  \item The homology groups  $H_m(\widetilde{M}(n))$ are free abelian. For the Betti numbers we have:
  $$\beta^m(\widetilde{M}(n))= \beta^m({M}(n))+\beta^{m-2}({M}(n)).\qed$$
\end{enumerate}

\end{corollary}

\begin{definition} \label{Dfn_areaR3}  Let  $(P,\xi) \in \widetilde{M}_3(n)$, let  $ (x_i,
y_i)$ be the vertices of $P$. The \textit{ area} of the pair $(P,\xi)$ is defined as
the scalar product:
$$2S(P,\xi)=( p_1 \times p_2 + p_2 \times p_3+ \dots + p_n\times p_1,\xi).$$

An alternative equivalent definition is:
$$S(P,\xi)=A(pr_{\xi^\perp}(P)),$$
where $pr_{\xi^\perp}$ is the plane orthogonal to $\xi$ and
cooriented by $\xi$.

\end{definition}

\medskip

\begin{proposition}\label{Thm_crirical_3D}
For an equilateral polygon with odd number of edges, critical points $(P, \xi)$ of the function $S$ are
 pairs $(P, \xi)$ such that $P$ is a planar cyclic
    polygon, and $\xi$ is orthogonal to the affine hull of $P$.
   If $(P, \xi)$ is a critical point, then $(P, -\xi)$ is critical as well.

\end{proposition}

Proof.
The paper  \cite{KhristofPanina} contains a characterization of all critical points for a generic polygonal linkage (which is not necessarily equilateral). In
our particular case it implies that
critical points $(P, \xi)$ of the function $S$ fall into two classes:
 \begin{enumerate}
    \item \textbf{Planar cyclic configurations.}
    These are pairs $(P, \xi)$ such that $P$ is a planar cyclic
    polygon, and $\xi$ is orthogonal to the affine hull of $P$.
    \item \textbf{Non-planar configurations.} They are characterized
    by the three following conditions:
    \begin{enumerate}
\item The vectors $\xi$ and $\overrightarrow{S}=p_1 \times p_2 + p_2 \times p_3+ \dots + p_n\times p_1$
are parallel (but they can have opposite directions).
    \item The orthogonal projection of $P$ onto the plane $\overrightarrow{S(P)}^\perp$ is a
    cyclic polygon.
    \item For every $i$, the vectors $\overrightarrow{T_i}$, $\overrightarrow{S}$, and $\overrightarrow{d_i}$ are
    coplanar.
\end{enumerate}
Here $\overrightarrow{d_i}$ is the $i$-th short diagonal,
$\overrightarrow{T_i}$ is the vector area of the triangle
$p_{i-1}p_ip_{i+1}$.

 \end{enumerate}

Let us show that the second class (non-planar configurations)  is empty.
Indeed, given  a non-planar critical configuration, introduce a Cartezian system with the $z$-axes parallel to $\xi$. The three conditions (a), (b), and (c) implie that the absolute value of the slope
of an edge with respect to the plane $(x,y)= \xi^\perp$ does not depend
on the edge. This implies a contradiction with the closing condition: $\sum_{i=o}^n (z(p_i)-z(p_{i-1}))=0,$ where the indices  are modulo $n$.\qed

\bigskip

\begin{theorem}\label{Cor_perfectMorse}
\begin{enumerate}
  \item The function $S$ is a perfect Morse function on the decorated configuration space $\widetilde{M}_3(L)$ for an equilateral linkage with odd number of edges.
  \item  Each  critical point of the function $S$ is a pair $(P,\xi)$, where
$P$ is a planar cyclic configuration, $\xi$ is a unit vector orthogonal to $P$.
Each planar cyclic configuration  $P$ gives two critical points of the function $S$  (with two different choices of the normal vector $\xi$).
\item The Morse index of a critical point $(P,\xi)$ is
$$m(P,\xi)=2e-2\omega -2,$$
where $\omega$ is the  winding number $P$ around the center of the circumscribed circle,$e$ is the number of edges that go counterclockwise.\footnote{ The vector $\xi$ sets an orientation on the plane of the polygon,
so it makes sense to speak of "edges going clockwise" and "edges going counterclockwise". }
\end{enumerate}

\end{theorem}
Proof. (\textit{i})
The second statement is already proven. Let us show that the number of critical points equals the sum of Betti numbers.
The latter
are already known due to A. Klyachko \cite{klya}:

$$\beta^{2p}(M_3(n)) = \sum _{0\leq i\leq
p}\left(\!\! \begin{array}{c}2k\\ i
\end{array}\!\!\right),\;\;p<k.$$

By Corollary \ref{Cor_Gisin}
 we
have
$$\beta^{2p}(\widetilde{M}_3(n))=\sum _{0\leq i\leq p}\left(\!\!
\begin{array}{c}n\\ i
\end{array}\!\!\right),\;\;p<k.$$

Each equilateral cyclic $n$-gon with an orthogonal vector $\xi$ is defined by its winding number $\omega$ and by the set of edges that go clockwise.
Assume that the winding number is positive (negative values are treated by symmetry).
If the number of edges that go clockwise is $e$, then the winding number  ranges from $  1$ to $ (k-e)$.

 For $p=0,1,..., k$ denote by $N^p_n$ the number of such cyclic
equilateral polygons for which
$e-\omega-1=p.$
Then
$\widetilde{\beta}^{2p}_n=N^p_n.$

(\textit{iii}) Straightforward analysis of the Hesian matrix is very complicated (probably impossible), so we use derive a combinatorial approach. We prove the formula for Morse index by induction by $n$. The base is given by equilateral pentagon.
By symmetry, we  assume in whatfollows that we have a critical point $(P,\xi)$ such that $S(P,\xi)>0$, or, equivalently, with winding number $\omega >0$.

Proof is based on two observations:

(1) There is a natural embedding of a neighborhood of $P$ in  space $M_2(n)$ to $\widetilde{M}_3(n)$. It sends a configuration $P$ to $(P,\xi)$ with the same $P$ and $\xi$ orthogonal to the affine hull of $P$. The direction of $\xi$ we choose in the way such that $S$ is positive. Let us choose a basis of the tangent space $T_{(P,\xi)}\widetilde{M}_3(n)$ which starts by the pushforward of a basis of $T_PM_2(n)$ and ends by coordinates of $\xi$. The Hessian matrix related to this basis is a block matrix:
$$HESS(P,\xi)=\left(
    \begin{array}{ccc}
      H_1 & 0  & 0\\
      0 & H_2 & 0 \\
      0 & 0  & -E\\
    \end{array}
  \right),
$$
where $H_1$ is the $(n-3)\times (n-3)$ Hessian matrix of the planar polygon $P$ related to the area function $A$ and the space $M_2(n)$; $E$ is the the $2\times 2$ unit matrix.

(2) For an $(n+2)$-gon and a number $1 \leq i \leq n$, consider the embedding $\varphi_i:\widetilde{M}_3(n)\rightarrow \widetilde{M}_3(n+2)$ which keeps $\xi$ and replaces the edge number $i$ by a fold of three edges (see Fig. 1, (3)  for an equilateral triangle with an edge by a three-fold).  Each critical point $(P,\xi)\in {\widetilde{M}_3}(n)$ induces a critical point $\varphi_i(P,\xi) \in \widetilde{M}_3(n+2)$. Since the embedding has codimension two, and all Morse indices can only be even, we have either
$m(\varphi_i(P,\xi)) = m(P,\xi)$, $m(\varphi_i(P,\xi)) = 2+ m(P,\xi)$, or $m(\varphi_i(P,\xi)) = 4+ m(P,\xi)$.
More precisely, replacement of an edge replaced by a three-fold adds two extra columns to $H_1$ and two extra columns to $H_2$.
We know from \cite{khipan} that the Morse index of $P$ related to $M_2(n)$ and the oriented area $A$ (that is, the number of negative eigenvalues of $H_1$) equals $e-2\omega-1$. It  increases by one after replacement of an edge by a three-fold. This only leaves the case  $m(\varphi_i(P,\xi)) = 2+ m(P,\xi)$.

These arguments allows to make an induction step $n \rightarrow n+2$ and thus proves (\textit{iii})  for all cyclic configurations with triple edges.

It remains to prove the formula for  configurations  without triple edges, that is, with all the edges going counterclockwise. We keep assuming that $S(P,\xi)>0$, so there are exactly $k$ such polygons: with $\omega =1,2,3,...,k$.
Their Morse indices should be $2n-4$, $2n-6$, $2n-8$, etc. The only question is which configuration has one or the other index.
We know that $H_1$ contributes $n-3$, $n-5$, $n-7$, etc. to each of the Morse indices. The block $H_2$ contributes not more than $n-3$,  and $-E$ contributes exactly two. The statement follows.
\qed

As an illustration, we list below
 all cyclic equilateral pentagons and heptagons. The first column depicts a combinatorial type, the third one
tells the number of configurations of this type, and the last one tells the Morse index.

\begin{figure}[h]
\centering
 \subfigure[Critical equilateral pentagons]{
 \includegraphics[width=6 cm]{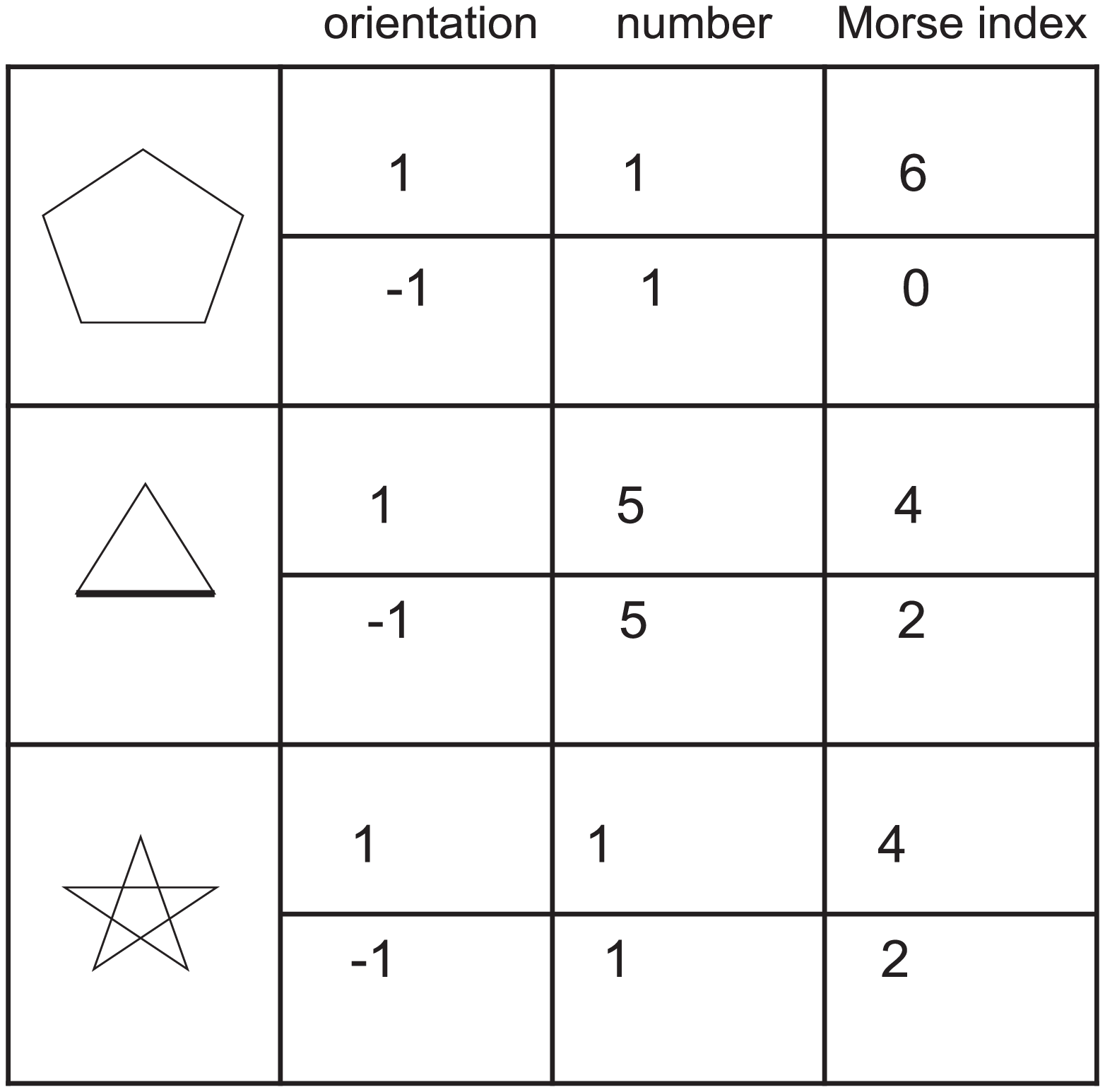}}
   \subfigure[Critical equilateral heptagons]{
    \includegraphics[width=6 cm]{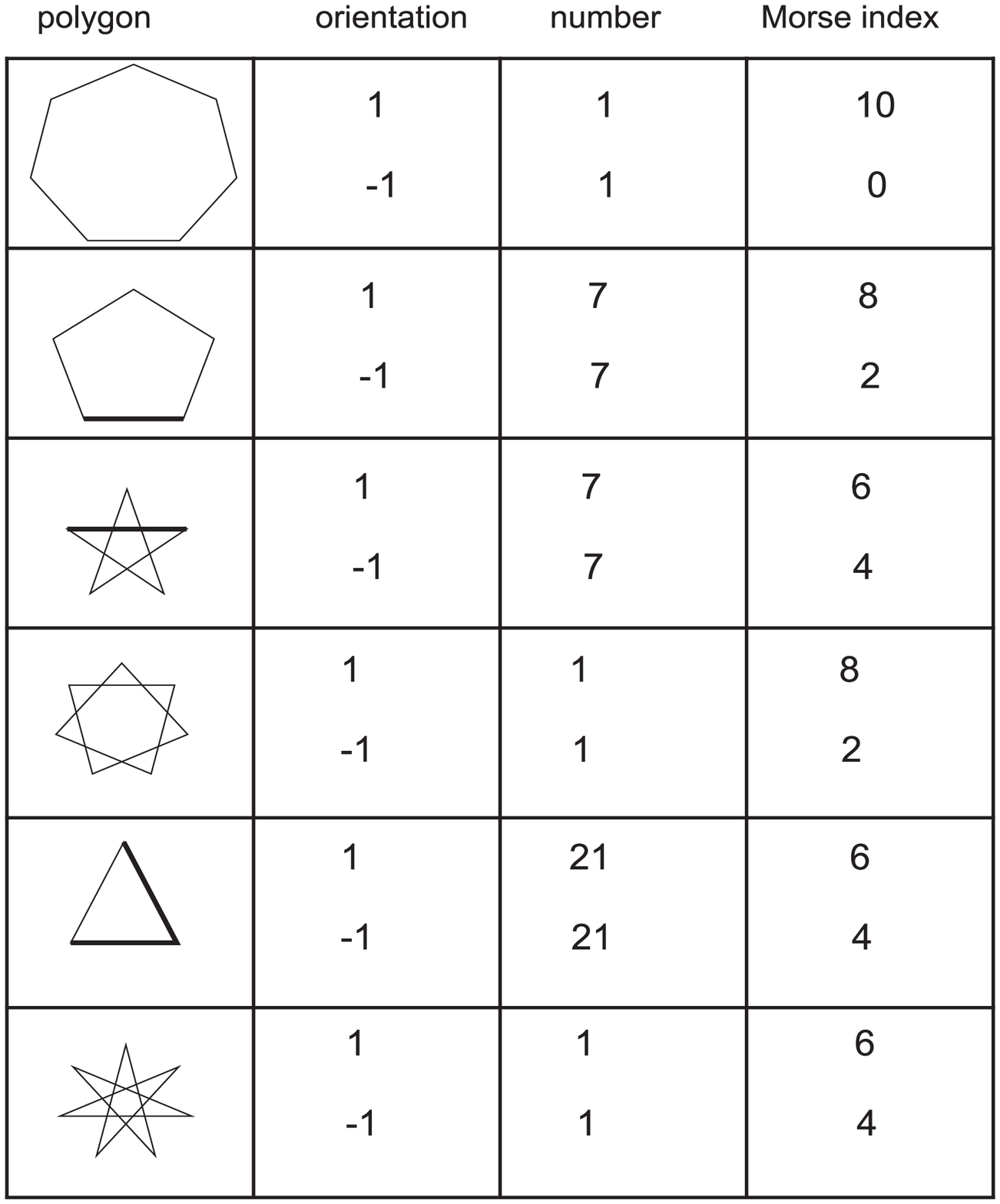}}

\label{Figure_7gon}
\end{figure}

\section*{Concluding remarks}
 Decorated configuration space and the ares $S$ can be defined for a polygonal linkage that is not necessarily equilateral.
  However, generically, the function $S$ is not a perfect Morse function.

In the paper we omit the discussion about non-degenericity of critical  points, since it appears to be somewhat technical. However the following arguments work:
We may  replace the equilateral linkage $(1,1,...,1)$ by its perturbation $(1+\epsilon_1,1+\epsilon_2,...,1+\epsilon_n)$  has non-degenerate critical configurations that are close to above described equilateral ones.

\section*{Acknowledgements} This work is
supported by the Russian Science Foundation under grant  16-11-10039.

\end{document}